\documentclass{article}
\usepackage{a4wide}
\usepackage{amssymb}
\usepackage{amsfonts}
\usepackage{amsxtra}
\usepackage{latexsym}

\begin{document}

\newtheorem{theo}{Theorem}[section]
\newtheorem{theo*}{Theorem}
\newtheorem{lemma}[theo]{Lemma}
\newtheorem{definition}[theo]{Definition}
\newtheorem{remark}[theo]{Remark}
\newtheorem{cor}[theo]{Corollary}
\newtheorem{prop}[theo]{Proposition}
\newtheorem{example}[theo]{Example}

\newcommand{\N}{{\mathbb{N}}}
\newcommand{\Z}{{\mathbb{Z}}}   
\newcommand{\Q}{{\mathbb{Q}}}
\newcommand{\R}{{\mathbb{R}}}
\newcommand{\C}{{\mathbb{C}}}
\newcommand{\M}{{\mathbb{M}}}
\newcommand{\F}{{\mathbb{F}}}
\newcommand{\LL}{{\mathbb{L}}}

\newcommand{\g}[1]{{\mathfrak {#1}}}
\newcommand{\cc}[1]{{\cal {#1}}}
\newcommand{\ul}{\underline }
\newcommand{\ol}{\overline }

\newcommand{\qed}{\hspace*{\fill}$\Box$}
\newcommand{\Qed}{\hspace*{\fill}$\Box \Box$}

\newcommand{\md}{{\rm mod}\,\,}
\newcommand{\ee}{{\rm zip}\,}
\newcommand{\ei}{{\rm zap}\,}


\title{Dynamical systems with double recursion are undecidable}
\author{Mihai Prunescu \thanks{ Brain Products, Freiburg, Germany, and Institute of Mathematics of the 
Romanian Academy, Bucharest, Romania.
mihai.prunescu@math.uni-freiburg.de.}}
\date{}
\maketitle

\begin{abstract}

\parindent 0 cm 

A primitive type of dynamical system is introduced. It is shown that there is no decision procedure able to answer if any
such dynamical system is ultimately zero.

\thanks{A.M.S.-Classification: 12L05.}
\end{abstract} 


\section{Introduction}

\begin{definition}\rm A {\bf dynamical system with double recursion} is a tuple $\g A = (A, f, 0, 1)$ 
consisting of: a finite set $A$, 
a function $f : A \times A \rightarrow A$, an element $1 \in A$ called start symbol and an element $0 \in A$ called white
colour. By {\bf system} or {\bf dynamical system} we shall understand a dynamical system with double recursion, if not 
otherwise is explicitly stated. If a dynamical system $\g S \models \forall \, x,y \,\, f(x,y) = f(y,x)$ 
we say that the system $\g S$ is {\bf symmetric}.

\end{definition}

\begin{definition}\rm The development of $\g A$ is a function $a : \N \times \N \rightarrow A$ defined as follows:
\[ a(i,j) = 
\begin{cases}
1 & {\rm if} \,\, i = 0 \,\, \vee \,\,  j = 0 ,\\
f(a(i-1,j),a(i,j-1)) & {\rm if} \,\, i > 0 \,\, \wedge \,\, j > 0.
\end{cases}
\]
If $\g A = \g S$ is symmetric, the development is always a symmetric matrix: $a(i,j) = a(j,i)$.
\end{definition}

\begin{definition}\rm
A dynamical system with double recursion $\g A$ is called {\bf ultimately zero} if 
\[ \mathfrak A \,\,\models \,\,\exists \, N \in \N \,\,\forall \, i,j \in \N \,\,\,\, i > 0 \, \wedge \, j > 0 \,
\wedge \, i+j > N \,\, \longrightarrow \,\, a(i,j) = 0.\]
The problem UW: given a dynamical system $\g A$, tell if it is ultimately zero.

The problem UWS: given a symmetric dynamical system $\g S$, tell if it is ultimately zero.

\end{definition}

\begin{theo} \label{main}

The problems {\rm UW} and {\rm SUW} are not algorithmically solvable.

\end{theo}

Of course, if SUW is not algorithmically solvable then neither shall be UW. However, I believe that starting with an
easier sketch of proof for UW is helpful for the natural introduction in the matter. 

For encoding the evolution of a Turing machine in the development of a dynamical system, 
we shall need the following description of the development:

\begin{definition}

The development $a: \N \times \N \rightarrow A$ of $\g A$ shall be seen as an infinite matrix with $a(0,0) = a_{0,0}$
in the left upper corner. Denote by:
\[
D_n = \{ a_{i,j} \,|\, i+j=n\}
\]
the $n$-th diagonal. Then $a = \cup D_n$ is a partition of the development of $\g A$.
\end{definition}

Before starting with the proof, we observe that dynamical systems can be completely algebraically modeled over finite
fields. Just take some prime-power $q > |\,A\,|$, fix an embedding $A \subseteq \F_q$ such that 
$0$ and $1$ are now the corresponding constants in $\F_q$, and find some polynomial $F \in \F_q[x,y]$ 
whose associated polynomial function is an extension of $f$. For $f$ symmetric one can always find a symmetric polynomial 
$F$. In this way the Theorem \ref{main} gets a number-theoretic meaning:

\begin{cor}\label{Numbertheory}
Consider dynamical systems of the form $\g F = (\F_q, F, 0, 1)$ where $\F_q$ are  
finite fields and $F \in \F_q[x,y]$ are polynomials. 
Then it is undecidable if dynamical systems $\g F$ are ultimately zero.
Moreover, this question remains undecidable if simultaneously restricted 
to the prime fields $\F_p$ and to symmetric polynomials.
\end{cor}

{\bf Proof}: For a proof one should recall at most the interpolation with more than one  variables. 
Take some arbitrary function
$f : \F_q^n \rightarrow \F_q$. To show that $f$ is always polynomial, we will write it in the form:
\[ f(\vec x) = \sum \limits _{\vec a \in \F_q^n} d_{\vec a} (\vec x) f(\vec a),\]
where $d_{vec a}(\vec x) = 0$ if $\vec x \neq \vec a$ and $d_{\vec a}(\vec a) = 1$. For $n=1$ take:
\[ d_a(x) = \prod \limits _{b \in \F_q \setminus \{a\}}\frac{ (x-b)}{(a-b)}.\]
Take $d_{\vec a} (\vec x) = d_{a_1}(x_1)d_{a_2}(x_2) \dots d_{a_n}(x_n)$. \qed

Some historical words. One can see a dynamical system with double recursion as the tiling of the quarter of plane with tiles
obtained by modifying squares such that the only relevant compatibilities to solve for puting down a 
new stone are the neighbors from North and West. Seen in this way, these results are related with the classical results
concerning undecidability of questions about finite sets of tiles, see \cite{BGG}, \cite{HW}, \cite{RR}. Aperiodic tilings
related to some linear dynamical systems over finite fields have been described by the author in \cite{MP}, another 
recent preprint. 

The unique classical ingredient used here is the Theorem of Rice, see \cite{GR}, in its modern formulation concerning
sets of (codes of) Turing machines, as stated for example in \cite{US}.


\section{The problem UW}

In order to prove that the problem UW is not algorithmically solvable we will interpret instances of the Halting Problem
in instances of UW. 

\begin{definition}\rm
An instance of the Halting Problem is a pair $(M,w)$ where $M=(\Sigma, Q, q_0, q_s, \bar b, \delta)$ is a Turing machine and
$w \in \Sigma^*$ is an input for $M$. Here the tape of $M$ is infinite in both directions, $\Sigma$ is the alphabet of $M$,
$Q$ is $M$'s set of states, $q_0$ and $q_s$ are the start state and respectively the stop state, $\bar b \in \Sigma$ is
the blank symbol, and $\delta : \Sigma \times Q \rightarrow \Sigma \times Q \times \{R, L, S\}$ is the transition function.
\end{definition}

\begin{lemma}\label{simple}
To every instance $(M,w)$ of the Halting Problem one can algorithmically associate a dynamical system $\g A = (A,f,0,1)$
such that $\g A \in$ UW if and only if for input $w$ the machine $M$ stops and after stopping the tape contains only blanc
symbols.
\end{lemma}

{\bf Proof}: If the Lemma \ref{simple} is true then UW is not algorithimically solvable. This is true because 
to stop with a clear tape is an undecidable property, according to the Theorem of Rice.

Let $(M,w)$ be an instance of the Halting Problem. The dynamical system $\g A$ shall be constructed step by step following 
its development:

\begin{itemize}
\item The start symbol of $\g A$ is a new letter $1$ that has nothing to do with the Turing machine. On the other hand 
the colour white $0$ of $\g A$ and the blanc symbol $\bar b$ of $M$ are the same.
\item Let the input $w \in \Sigma ^*$ be the word $w_1 \dots w_n$.
Using a set $U$ of new letters one defines $f$ such that the content of a diagonal $D_w$ of the development is exactly:
\[
1 \,\,0\,\,0\,\,\delta_0 \,\, w_1 \, \, \dots \,\, w_n \,\, 0 \,\, 0 \,\, 1
\]
The letters of $U$ shall be used only for this goal, and then never again. 
The simulation of the Turing machine starts with this diagonal.
\item The just constructed diagonal $D_w$ is said to be a diagonal of type $0$. Starting with $D_w$ diagonals are
alternatively of types $0$, $1$, $0$, $1$, and so on. Successive diagonals of type $0$ simulates successive configurations
of the Turing machine. The diagonals of type $1$ between them are used to transfer the information.
\item The alphabet used for diagonals of type $0$ contains 
$\Sigma \cup (\Sigma \times Q)$. We denote letters from $\Sigma \times Q$
with $\delta_i$. The meaning of the letter "$(a,q)$" is that the head of $M$ reads $a$ with $M$ being in the state $q$.
\item If $\Gamma_0$ is the alphabet used for diagonals of type $0$, the alphabet used for diagonals of type $1$ will be
$\Gamma_1 =(\Gamma_0 \times \Gamma_0 \setminus \{(0,0)\}) \cup \{0\}$.
\item The function $f$ is defined on diagonals of type $0$ in the obvious way $f(a,b) = (a,b)$ if at least one of $a$ and
$b$ are not $0$ or $1$, and $f(0,0) = f(1,0) = f(0,1) = 0$.
\item The function $f$ is defined on diagonals of type $1$ such that: if the element $a_{i,j}$ of the last diagonal of type
$0$ simulated a certain cell of the tape of $M$ at a given time $k$, then the cell $a_{i+1,j+1}$ simulates the same cell
at the time $k+1$. The following example shows how the diagonal of type $1$ in between makes possible that
the element $a_{i+1,j+1}$ gets information from three successive cells: $a_{i,j}$ and the elements simulating their 
neighbors.
\[
\begin{array}{ccc}
    &    	&   b 	      \\
    &  \delta   & (\delta, b) \\
 a  &(a,\delta) &   c
\end{array}
\]
\item Every diagonal of type $0$ is with two cells longer then the precedent one
and the head makes one step per time, so there is no danger that the simulation
leaves the carpet or even that the simulation meats the wall of ones.
\item For the special letter $\delta = (0,q_s)$ we define $f$ such that $f(\delta,0)=f(0,\delta)=0$. This makes the
development ultimately zero if and only if the machine stops with clear tape.
\item Now take $A$ to be $\{1\} \cup U \cup \Gamma_0 \cup \Gamma_1$ and $f : A \times A \rightarrow A$ to respect all the
conditions given above.

\end{itemize} 

\qed


\section{The problem SUW}\label{SUW}

\begin{definition}\rm
Let $\Gamma \neq \emptyset$ be a set and $\equiv $ be the partition of $\Gamma \times \Gamma$ consisting of the following 
sets: for all $a \in \Gamma$ the singleton sets $\{(a,a)\}$ and for all $a, b \in \Gamma$ with $a \neq b$ the two-element
sets $\{(a,b), (b,a)\}$. Then $\equiv$ is an equivalence relation over $\Gamma \times \Gamma$. Consider the set
of equivalence classes:
\[ \Gamma \cdot \Gamma = (\Gamma \times \Gamma) / \equiv \]
which is the set of unordered pairs of elements of $\Gamma$. We denote the equivalence class of $(a,b)$ with $[a,b]$ and call
this the unordered pair of $a$ and $b$.
\end{definition}

\begin{lemma}\label{hard}
To every instance $(M,w)$ of the Halting Problem one can algorithmically associate a 
symmetric dynamical system $\g S = (S,f,0,1)$
such that $\g S \in$ SUW if and only if for input $w$: 
(the machine $M$ stops with white tape without having done any step in the negative side of the tape) or (the machine $M$
shall make a step in the negative side of its tape and the first time when $M$ makes such a step the tape of $M$ is clear).
\end{lemma}

{\bf Proof}: If the Lemma \ref{simple} is true then SUW is not algorithimically solvable. This is true because 
the given condition is an undecidable property, according to the Theorem of Rice.

Before starting the construction, I shall shortly explain the arrising difficulties. We construct again the function $f$
together with its development. The function $f$ and the development are both symmetric, so we will only consider the 
right half of the development. It shall be again so, that on some special (half of) diagonals one simulates successive
configurations of the Turing machine on input $w$. The function $f$ being symmetric, one cannot directly make the difference
between Left and Right. To overcome this difficulty one can try to double the number of letters of $\Gamma$ and write every
letter $c$ as $cc'$. The function $f$ should now act symmetrically on diagonals of type $0$, so we define $f(a,b)$ to be
the unordered pair $[a,b]$. This strategy is not sophisticated enough: if we look at words $aba$ and $bab$ on a diagonal
of type $0$, they both produce a word $xx$ on the following diagonal of type $1$, where $x = [a,b]$. This means that this
encoding may lose essential information. The solution shall be to triple the number of letters and to encode every letter
$c$ with a sequence $cc'c''$, where $c'$ and $c''$ are used only for this goal. 

\begin{itemize}
\item The start symbol of $\g S$ is a new letter $1$ that has nothing to do with the Turing machine. Also the 
the colour white $0$ of $\g S$ is now a new letter.
\item Let the input $w \in \Sigma ^*$ be the word $w_1 \dots w_n$.
Using a set $U$ of new letters one defines $f$ symmetrically 
such that the content of a diagonal $D_w$ of the development is exactly:
\[
1 0^9 \,\, w_n''w_n'w_n\,\,\dots \,\, w_1''w_1'w_1\,\,\delta_0 '' \delta_0 ' \delta_0 \,\,000\,\,
\delta_0\delta_0'\delta_0''\,\, w_1w_1'w_1'' \, \, \dots \,\, w_n w_n'w_n''\,\, 0^9 1
\]
Here means $0^9$ a word built up by $9$ zeros. The letters of $U$ shall be used only for this goal, and then never again. 
The simulation of the Turing machine starts with this diagonal.
\item The just constructed diagonal $D_w$ is again a diagonal of type $0$. This time there are $8$ types of
diagonals: types $0$, $1$, $\dots$, $7$.  Starting with $D_w$ diagonals are
of types $0$, $1$, $\dots$, $7$, $0$, $1$, $\dots$, $7$, and so on. 
Successive diagonals of type $0$ simulates successive configurations
of the Turing machine. The diagonals of types $1$, $\dots$, $7$ between them are used to transfer the information from
a simulation to the next.

\item The alphabet $\Gamma_0$ used for diagonals of type $0$ contains three disjoint copies of the set
$\Sigma \cup (\Sigma \times Q) \setminus \{\bar b\}$. The letter $c \in \Sigma \setminus \{\bar b\}
$ is called $c'$ in $\Sigma'$ and
$c''$ in $\Sigma''$, the same for the letters $\delta \in \Sigma \times Q$. 
The fact that $c \in \Sigma$ is contained in a cell of the Turing machine is encoded by the word $cc'c''$ on a diagonal
of type $0$. On the right side of the development: 
The fact that the head of the Turing machine in state $q$ reads a cell of content $a$ 
is encoded by the word $\delta \delta' \delta''$ on a diagonal of type $0$, where  $\delta = (a,q) \in \Sigma \times Q$ . 
The blanc symbol $\bar b$ as content of a cell of the Turing machine is always encoded in the simulation by the word $000$
on a diagonal of type $0$. In the left side of the development the codes are $c''c'c$, $\delta''\delta'\delta$ and $000$
respectively.
\item Let $\Gamma_0$ be the alphabet used for diagonals of type $0$. 
For $i = 1, 2, \dots, 7$ the alphabet used for diagonals of type $i$ will be
$\Gamma_i =(\Gamma_{i-1} \cdot \Gamma_{i-1} \setminus \{[0,0]\}) \cup \{0\}$.
\item The function $f$ is defined on diagonals of type $i = 0,1,\dots, 6$ in the obvious way $f(a,b) = [a,b]$ 
if at least one of $a$ and $b$ are not $0$ or $1$, and $f(0,0) = f(1,0) = f(0,1) = 0$.
\item The function $f$ is defined on diagonals of type $7$ such that: if the element $a_{i,j}$ of the last diagonal of type
$0$ contains a letter $c$, $c'$, $c''$, $\delta$, $\delta'$, $\delta''$ or $0$ that appears
in a subword of length $3$ simulating a cell of the tape of $M$ at a given time $k$, then the element $a_{i+4,j+4}$ 
of the development
shall be the corresponding letter of the diagonal word simulating the configuration of $M$ at time $k+1$.
This is done like in the following example. Let $a_{i,j} = \delta' \in \Gamma_0$ be a part of the following segment of 
simulation in a diagonal of type $0$: $\dots cc'c''\delta \delta' \delta'' dd'd''\dots$,
and suppose that in the next configuration the tape-cell containing $\delta$ shall contain $e \in \Sigma$. 
As proved in the postponed Lemmas \ref{symcod} one has on the first coming diagonal of type 
$7$: $a_{i+3,j+4} = \alpha$ and $a_{i+4,j+3} = \beta$ with
$\alpha, \beta \in \Gamma_7$, such that $\alpha$ encodes the word $cc'c''\delta\delta'\delta''dd'$ or its reverse and
$\beta$ encodes the word $c'c''\delta\delta'\delta''dd'd''$ or its reverse. One has either the words themselves 
(if we look to the right-hand side of the development), or the
reversed words (if we look to the left-hand side of the development). If we are in the right-hand side the matching of the 
encoded words looks like:
\[
\begin{array}{ccccccccc}
  c  &    c'	&   c''  & \delta & \delta' & \delta '' & d & d' &	      \\
     &     c'   &   c''  & \delta & \delta' & \delta''  & d & d' & d''
\end{array}
\]
If we are in the left-hand side of the development, the matching is:
\[
\begin{array}{ccccccccc}
  d''   &    d'	&   d  & \delta'' & \delta' & \delta  & c'' & c' &	      \\
        &    d' &   d  & \delta'' & \delta' & \delta  & c'' & c' & c
\end{array}
\]
In both cases the matching is a word of length $7$ centrated in $\delta'$, so the value $f(\alpha,\beta)$ is uniquely
determined to be $e'$, where $e \in \Sigma$ is the letter that shall replace
 $\delta$ in that tape-cell in the next configuration. 
The same arguments work for every connected subword of length $8$ that is disjoint from the central $000$ word,
like $c''\delta \delta'\delta''dd'd''e$, and so on.  
\item Every diagonal of type $0$ is with eight elements longer then its predecessor of type $0$ (four elements
in the left-hand side and four elements in the right-hand side)
and the simulation needs at most three  elements more per step, so there is no danger that the simulation
leaves the carpet or even that the simulation meats the wall of ones.

\item If the connected subword $\delta \delta'\delta''$ with the special letter $\delta = (0,q_s)$ arrises, then 
we define $f$ such that the corresponding development elements in the next diagonal of type $0$ are $0$.

\item For the connected subwords of length $8$ containing the central $000$
 the funcion $f$ is defined such that: Words of the type $c''c'c000cc'c''$ are preserved in the next
configuration. Words $\delta''\delta'\delta000\delta\delta'\delta''$ 
are replaced with $e''e'e000ee'e''$ if $\delta = (a,q) \rightarrow (e,q',R)$. 
Words $\delta''\delta'\delta000\delta\delta'\delta''$ are replaced with $a''a'a000aa'a''$ if
$\delta = (a,q) \rightarrow (e,q',L)$. As asked in the statement of this Lemma the computation dies 
by the first movement in the left-hand side of the tape.
 
\item Now take $A$ to be $\{1\} \cup U \cup   _{i=0} ^ 7 \Gamma_i$ and $f : A \times A \rightarrow A$ to respect all the
conditions given above.

\end{itemize} 

\qed


\section{Symmetric codes}

\begin{definition}\rm
Let $\Gamma_0$ be a finite alphabet with $\geq 2$ letters and $0\in \Gamma_0$ a special letter. 
We define the sequence of alphabets $\Gamma_i$ such that $\Gamma_{i+1} =( \Gamma_i \cdot \Gamma_i \setminus \{[0,0]\})
\cup \{0\}$ and $f:\Gamma_i \times \Gamma_i \rightarrow \Gamma_{i+1}$ such that $f(a,b)=[a,b]$ if at least one of
the arguments is not $0$ and $f(0,0)=0$. For an alphabet $\Gamma$ let $\Gamma^*$ be the set of words over $\Gamma$ and
$\Gamma^{\geq k}$ the set of words of length $\geq k$ over $\Gamma$. The set of words of length $k$ shall be simply denoted
$\Gamma^k$.
\end{definition}

\begin{definition}\rm
Let $\pi_i : \Gamma_i^{\geq 2} \rightarrow \Gamma_{i+1}^*$ given as $\pi_i(w_1\dots w_n)=f(w_1,w_2)f(w_2,w_3)\dots f(w_{n-1},
w_n)$. Let $\pi : \Gamma_0^8\rightarrow \Gamma_7$ given as $\pi(w) = \pi_6\pi_5\pi_4\pi_3\pi_2\pi_1\pi_0(w)$. We call
the words $w$, $\pi_0(w)$, $\pi_1 \pi_0(w)$, $\dots$, $\pi_5\pi_4\pi_3\pi_2\pi_1\pi_0(w)$ the levels of the code.
\end{definition}

\begin{definition}\rm
Now let $\Gamma_0$ be the alphabet defined in the section \ref{SUW}. Let $E$ be the set of all words in $\Gamma_0^8$
that can possible arrise during a simulation. They are exactly the connected subwords of length $8$ in all words 
$aa'a''bb'b''cc'c''dd'd''$ where $a,b,c,d \in \Sigma \cup \Sigma \times Q$ are not necessarily different, 
and if some $e \in \{a,b,c,d\}$ are $0$ then the corresponding $e'=e''=0$. The restriction of $\pi$ to 
$E \rightarrow \Gamma_7$ shall be simply called $\pi$.
\end{definition}

\begin{definition}\rm
Let $S$ be the set of connected subwords of length $8$ in all words $ba''a'a000aa'a''b$ where 
$a, b \in \Sigma \cup \Sigma \times Q$ are not necessarilly different, and if $a = 0$  then the corresponding $a'=a''=0$. 
Again the restriction of $\pi$ to 
$S \rightarrow \Gamma_7$ shall be simply called $\pi$.
\end{definition}

\begin{definition}\rm
For a word $w \in \Gamma^*$, $w=w_1\dots w_n$, call $\sigma(w)$ the reversed word $w_n\dots w_1$.
\end{definition}

\begin{lemma}\label{symcod}
For all words $v \in E \cup S$ and $w \in (\Sigma \cup (\Sigma \times Q))^8$, 
if $\pi(w) = \pi(v)$ then $w = v$ or $w = \sigma(v)$.
\end{lemma}

{\bf Proof}: The proof works as follows: We check for all types of words in $E \cup S$ that only they and their reverses
lead to their symmetrically iterated codes.

{\it Words in $E$:}

It is enough to check the worst cases with letter repetitions. Start with the word $cc'c''cc'c''cc'$.
The levels of encoding are as follows: $x_1x_2x_3x_1x_2x_3x_1x_2$, $y_1y_2y_3y_1y_2y_3y_1$, 
$z_1z_2z_3z_1z_2z_3$, 
$v_1v_2v_3v_1v_2$, $t_1t_2t_3t_1$, $u_1u_2u_3$, $s_1s_2$, $\alpha$, where $\alpha \in \Gamma_7$. Starting with $\alpha$,
one gets $\alpha = [s_1,s_2]$ so the two possibilities are $s_1s_2$ and $s_2s_1$. The first one leads directly to $w$, 
the other one directly to $\sigma(w)$; there are not other possibilities to reconstruct the word.

For a complete proof of the Lemma, one has to check the following worst cases: (a) All the connected subwords of length $8$
in $cc'c''000cc'c''000000$ and (b) All the connected subwords of length $8$ in $cc'c''cc'c''cc'c''000$. 
All this cases have the following common property: at all levels
of coding, including the level $0$, two successive letters are equal if and only if they are $0$.

{\it Words in $S$:}

The words occuring here are the exceptions in our symmetric encoding: They do not enjoy the property that 
successive letters at every level are equal if and only if they are both $0$, but they are however well behaving even by
being the only words that don't enjoy this property. Look at $c''c'c000cc'$. The levels of its symmetric code are:
$x_1x_2x_300x_3x_2$, $y_1y_2y_30y_3y_2$, $z_1z_2z_3z_3z_2$, $v_1v_2v_3v_2$, $t_1t_2t_2$, $s_1s_2$, $\alpha$. In decoding
we have again the choice $s_1s_2$ or $s_2s_1$. If we choose $s_1s_2$, that can backwards 
develop only in $t_1t_2t_2$, and so on. One easily checks all other words in question.
\qed 

Just some commentaries at the end. Letters have been encoded by directed words of the form $cc'c''$ or 
$\delta\delta'\delta''$ in order to make the difference between the left and the right neighbor in a symmetric dynamical
system. To symmetrically encode words of length $8$ one needs $7$ supplementary types of diagonals. A letter of type $7$
encodes a word of length $8$ and so has information from three successive simulated Turing cells. The common part of two 
successive subwords of length $8$ has length $7$ and so always has a central letter: this is the letter to copy or replace
on the next coming diagonal of type $0$.

Both constructions done here work in polynomial time.

One can do very spectacular experiments with dynamical systems in three variables, developed with the rules 
$a_{i,0}= a_{0,i}=1$ and $a_{i,j} = f(a_{i-1,j},a_{i-1,j-1}, a_{i,j-1})$ and using for example functions $f(x,y,z)$ which
are symmetric in $x$ and $z$. This class of dynamical systems is also undecidable because it trivially contains the problem
SUW given above.

\end{document}